\documentclass[11pt]{article}
\usepackage[english]{babel}
\usepackage{times}
\usepackage{paralist}
\usepackage{amsmath,amstext,amssymb,amsthm, amsfonts}
\usepackage[dvips]{graphics,color}
\newtheorem{lemma}{Lemma}[section]
\newtheorem{theorem}[lemma]{Theorem}
\newtheorem{corollary}[lemma]{Corollary}

\newtheorem{proposition}[lemma]{Proposition}
\newtheorem{example}[lemma]{Example}

\def\X{\mathbb{X}}

\def\K{\mathbb{K}}

\def\Y{\mathbb{Y}}
\def\S{\mathbb{S}}
\def\R{\mathbb{R}}

\def\A{\mathbb{A}}
\def\M{\mathbb{M}}
\def\RR{\overline{\mathbb{R}}}
\def\eb{{\varepsilon_\beta}}
\def\gbm{{\gamma_{\beta,M}}}
\def\nbm{{n^*_{\beta,M}}}
\def\dbm{{\delta_{\beta,M}}}
\def\pbm{g_{\beta,M}}

\newcommand{\slim} {\mathop{\rm lim\,sup}}
\newcommand{\ilim} {\mathop{\rm lim\,inf}}

\title{MDPs with Setwise Continuous Transition Probabilities}
\begin{document}

\maketitle

\begin{center}
Eugene~A.~Feinberg \footnote{Department of Applied Mathematics and
Statistics,
 Stony Brook University,
Stony Brook, NY 11794-3600, USA, eugene.feinberg@sunysb.edu}, and
Pavlo~O.~Kasyanov\footnote{Institute for Applied System Analysis,
National Technical University of Ukraine ``Kyiv Polytechnic
Institute'', Peremogy ave., 37, build, 35, 03056, Kyiv, Ukraine,\
kasyanov@i.ua.}\\

\bigskip
\end{center}

\begin{abstract}
This paper describes the structure of optimal policies for infinite-state Markov Decision Processes with setwise continuous transition probabilities. The action sets may be noncompact. The objective criteria are either the expected total discounted and undiscounted costs or average costs per unit time.  The analysis of optimality equations and inequalities is based on the optimal selection theorem for inf-compact functions introduced in this paper.
\end{abstract}

\section{Introduction}
This paper studies infinite-state Markov Decision Processes (MDPs) with setwise continuous transition probabilities.  In order to ensure the existence of optimal policies and relevant properties for MDPs, such as validity of optimality equations and inequalities and convergence of value iterations, some continuity assumptions on transition probabilities and costs are required.  The two classic assumptions on transition probabilities, weak and setwise continuity of transition probabilities, were introduced for MDPs with compact action sets in \cite{Sc} for problems with expected discounted costs and used in \cite{Schal} for problems with average costs per unit time, where the additional Assumption~${\rm B}$ on the finiteness of relevant values functions was introduced.  The results from \cite{Schal} on average costs with setwise continuous transition probabilities were extended in \cite{HL} to problems with noncompact action sets. More general Assumption~${\rm \underline{B}}$ was introduced  in \cite{FKZ2012}, where the theory for average-costs criteria was developed for MDPs with weakly continuous transition probabilities and possibly noncompact action sets.  This paper provides results on average-cost MDPs with noncompact action sets and with setwise continuous transition probabilities satisfying Assumption~${\rm \underline{B}}.$ MDPs with noncompact action sets are important for several applications including inventory control \cite{Ftut} and linear quadratic stochastic control~\cite{FKZ2012, HL}.

Weak continuity of probabilities is a more general property than setwise continuity, and in some applications, including inventory control \cite[Section 4]{FLe07} and problems with incomplete information \cite{POMDP}, weak continuity of transition probabilities leads to results that cannot be achieved by applying models with setwise continuous transition probabilities.  However, models with weakly continuous probabilities require joint continuity properties of transition probabilities and costs in state-action pairs, while models with setwise continuous transition probabilities require continuity properties of transition probabilities and cost function only in the action parameter.  Because of this, models with weakly continuous transition probabilities deal only with problems with semi-continuous value functions, while models with setwise continuous transition probabilities can deal with problems with arbitrary measurable value functions.  For example, if action sets are finite, then models with setwise continuous transition probabilities are more general, and they were used in \cite{FP} to study MDPs with arbitrary measurable transition probabilities.

Optimality operators and equations play the central role for MDPs, and the analysis of MDPs with weakly continuous transition probabilities and possibly noncompact action sets  became possible after a classic fact in optimization, the Berge maximum theorem, was extended in \cite{FKZ2012, FKZ1} to noncompact action sets by introducing the notion of $\K$-inf-compact functions.  Berge's maximum theorem states continuity properties of optimality operators and equations.  For problems with setwise continuous transition probabilities, optimal selection theorems, that imply measurability properties of values, play the similar role; see \cite[p. 183]{HLL}.

In this paper we introduce an optimal  selection theorem for an inf-compact function (Theorem~\ref{th2}),  describe the theory for expected total discounted  and undiscounted costs in Theorem~\ref{th:dcoe},  establish in Theorem~\ref{teor1} the validity of optimality inequalities for average-cost MDPs satisfying Assumption~${\rm \underline{B}},$
and show in Section~\ref{s:exa} that Assumption~${\rm \underline{B}}$ is more general than Assumption ${\rm B}.$

\section{Optimal Selection Theorem}\label{s2}
Selection and optimal selection theorems play important roles in dynamic programming and in the theory of Markov Decision Processes (MDPs); see e.g., \cite{BS, DY, Ev, HLL, HPV}.  Selection theorems provide sufficient conditions for a graph of a set-valued function to contain a measurable function, called a selector, defined on the domain of the graph of the set-valued function.  Optimal selection theorems provide sufficient conditions for the selector to be an optimal solution for a parametric optimization problem.  In this section we provide an optimal selection theorem useful for the analysis of MDPs with setwise-continuous transition probabilities.

Let $\X$ and $\A$ be Borel spaces, that is, they are measurable subsets of Polish (complete, separable, metric) spaces, and let $\mathcal{B}(\X)$ and $\mathcal{B}(\A)$ be their Borel $\sigma$-algebras. A set-valued map $A:\X\to 2^\A$ is called strict if $A(x)\ne\emptyset$ for each $x\in\X.$ For a strict set-valued map $A:\X\to 2^\A,$ we define its graph  ${\rm Gr}_X(A):=\{(x,a)\in X\times\A: a\in A(x)\}$ restricted to $X\subset \X. $ When $X=\X$ we write ${\rm Gr}(A)$ instead of ${\rm Gr}_\X(A).$ A Borel mapping $\varphi:\X\to\A$ is called a measurable selector for  $A$ if $\varphi(x)\in A(x)$ for all $x\in\X.$

Let $\R$ denote the set of real numbers, $\R_+:=[0,+\infty),$ and $\RR=\R\cup\{+\infty\}$. We recall that a function $f:E\to\RR$ is called inf-compact on $E,$ where $E$ is a subset of a metric space, if the set $\mathcal{D}(\lambda,f):=\{e\in E:f(x)\le \lambda\}$ is compact for every $\lambda\in\R.$ A function $f:E\to\RR$ is called lower semi-continuous (l.s.c.) at $e\in E,$ if $f(e)\le \liminf_{n\to\infty} f(e_n)$ for all $e_n\to e$ with $e_n\in E$ for all $n\ge 1.$ A function $f$ is l.s.c. on $E$ if it is l.s.c. at all $e\in E.$ If a function $f$ is inf-compact on $E,$ then it is l.s.c. on $E.$ The function $f$ is upper semi-continuous (u.s.c.) at $e\in E$ (on $E$) 
if $-f$ is l.s.c. at $e\in E$ (on $E$). 
We denote by ${\rm dom}(f) :=\{e\in E\,:\, f(e)<+\infty\}$  the domain of $f{:\X\to\mathbb{R}\cup\{\pm\infty\}.}$ 

With a strict set-valued mapping $A:\X\to 2^\A$ and a function $u:{\rm Gr}(A)\to \RR,$ we associate the value function 
\begin{equation}\label{eq:val_func}
v(x):=\inf_{a\in A(x)} u(x,a),\quad\mbox{for all } x\in\X.
\end{equation}

The following theorems are useful for the analysis of MDPs with setwise continuous transition probabilities. We recall that a subset of a metric space is called $\sigma$-compact if it is a countable union of compact sets.
\begin{theorem}\label{th1}
Let $\X$ and $\A$ be Borel spaces, let $A:\X\to 2^\A$ be a strict set-valued map such that ${\rm Gr}(A)$ is Borel, and let $u:{\rm Gr}(A)\to \RR$ be a  Borel measurable function such that $u(x,\,\cdot\,)$ is l.s.c. on $A(x)$ for each $x\in\X.$ If each set $A^f(x)=\{a\in A(x): u(x,a)<+\infty\}$ is $\sigma$-compact, $x\in\X,$ then the value function $v$ defined in \eqref{eq:val_func} is Borel measurable.
\end{theorem}
\begin{proof}
For any set $Z\subset\X\times\A,$ let ${\rm proj}_\X(Z):= \{x\in\X\,:\,(x,y)\in Z \mbox{ for some }y\in \Y\}$ be the projection of $Z$ on $\X.$ Then ${\rm proj}_\X({\rm Gr}(A^f))= {\rm dom}(v).$  Since the function $u$ is measurable, the set ${\rm Gr}(A^f)={\rm dom}(u)$ is a measurable subset of $\X\times\A.$ Since each set $A^f(x)$ is $\sigma$-compact, the Arsenin-Kunugui theorem~\cite[Theorem 18.18]{Ke} or \cite[Theorem 1]{BP} implies that the set ${\rm dom}(v)$ is measurable.  According to~\cite[Corollary 1]{BP}, which follows from the Arsenin-Kunugui theorem, for every $n\ge 1$ there exists a measurable selector $\varphi_n:{\rm dom}(v) \to \A$ for $A^f$ such that $u(x,\varphi_n(x))\le v(x)+\frac{1}{n},$ if $-\infty<v(x)<+\infty,$  and  $u(x,\varphi_n(x))\le -n$  if $v(x)=-\infty.$  Therefore, $v(x)=\lim_{n\to\infty} u(x,\varphi_n(x))$ is a measurable function on ${\rm dom}(v).$ In addition, $v(x)=+\infty$ if $x\in \X\setminus {\rm dom}(v).$  Thus, the function $v:\X\to\RR$ is Borel measurable.
\end{proof}
\begin{theorem}\label{th2} {\rm (Optimal selection).}
Let $\X$ and $\A$ be Borel  spaces, let $A:\X\to 2^\A$ be a strict set-valued map such that ${\rm Gr}(A)$ is Borel, and let $u:{\rm Gr}(A)\to \RR$ be a  Borel measurable function such that $u(x,\,\cdot\,)$ is  inf-compact on $A(x)$ for each $x\in\X.$ Then ${\rm dom}(v)\in\mathcal{B}(\X),$ and there exists a Borel measurable selector $\varphi:{\rm dom}(v)\to\A$ for $A$ such that $u(x,\varphi(x))=v(x)$ for all $x\in{\rm dom}(v).$
Moreover, the value function $v$ {defined in \eqref{eq:val_func}} is Borel measurable.
\end{theorem}
\begin{proof}
The inf-compactness property of $u$ implies that $v(x)>-\infty$ for all $x\in\X,$  and $A^*(x):=\{a\in A(x): u(x,a)=v(x)\}$ is a compact set for each $x\in {\rm dom}(v).$ In view of Theorem~\ref{th1}, the function $v$ is Borel measurable and ${\rm dom} (v)\in\mathcal{B}(\X).$ Therefore, ${\rm Gr}(A^*)=\{(x,a)\in {\rm dom}(v)\times\A\,:\,  u(x,a)=v(x)\}\in {\cal B}(\X\times\A).$ In view of the Arsenin-Kunugui theorem, there is a measurable selector $\varphi:{\rm dom}(v)\to \A$ for $A^*.$  That is,  $u(x,\varphi(x))=v(x)$ for all $x\in{\rm dom}(v).$
\end{proof}
Theorem~\ref{th2} implies the following corollary.
\begin{corollary}\label{corsel:2}  {\rm (Optimal selection).}
Let assumptions of Theorem~\ref{th2} hold. If, additionally, one of the following conditions holds:
\begin{itemize}
\item[{\rm(i)}] there exists a Borel measurable selector $\varphi:\X\setminus {\rm dom}(v)\to\A$ for $A;$
\item[{\rm(ii)}] the set $\X\setminus {\rm dom}(v)$ is finite or countable;
\item[{\rm(iii)}]  $v(x)<+\infty$ for all $x\in\X,$ that is, for each $x\in \X$ there exists $a\in A(x)$ with $c(x,a)\in\R;$
\item[{\rm(iv)}] the function $u$ is real-valued;
\end{itemize}
then there exists a Borel measurable selector $\phi:\X\to\A$ for $A$ such that $u(x,\varphi(x))=v(x)$ for all $x\in \X.$
\end{corollary}
\begin{proof}
(i) According to Theorem~\ref{th1}, ${\rm dom}(v)\in \mathcal{B}(\X),$  and the measurable selector $\phi$ is defined at $x\in {\rm dom}(v).$ Set $\phi(x)=\varphi(x)$ for
$x\in \X\setminus {\rm dom}(v).$ (ii) follows from (i) since every selector $\varphi:\X\setminus {\rm dom}(v)\to\A$ for $A$ is measurable since the set $\X\setminus {\rm dom}(v)$   is finite or countable.  (iii) and (iv) follow from (ii) since $\X\setminus {\rm dom}(v)=\emptyset.$
\end{proof}
An l.s.c. function defined on a compact set is inf-compact.  Therefore, Corollary~\ref{corsel:2}(iv) implies \cite[Proposition D.5(a)]{HLL}, which was originally proved in \cite[Theorem 2]{HPV}, where the inf-compactness of $u(x,\cdot)$ on $A(x)$ is replaced with the assumptions that $u(x,\cdot)$ is l.s.c. on $A(x)$, and $A(x)$ are compact, $x\in\X.$ Corollary~\ref{corsel:2}(iv) also implies \cite[Proposition D.6(a)]{HLL}, where the additional assumption that the function $u$ is l.s.c. is imposed; this proposition is derived in \cite{HLL} from the results in \cite{Ri}.

\section{MDPs with Setwise Continuous Transition Probabilities}

Consider a discrete-time MDP specified by a tuple $(\X,\A,\{A(x):x\in\X\},c,q)$ with a Borel state space $\X,$ a Borel action space $\mathbb{A},$ nonempty Borel sets of feasible actions $A(x)$ of $\mathbb{A}$ at $x\in\X,$ one-step costs $c$, and
transition probabilities $q;$  see e.g., \cite{FKZ2012, HL, HLL, Schal}. The MDP satisfies the following standard assumptions: \begin{itemize}
\item[(a)] the graph of $A$ is measurable, that is, ${\rm Gr}(A)\in\mathcal{B}(\X\times\A);$
\item[(b)] there exists a measurable selector for $A:\X\to 2^\A;$
\item[(c)] the cost function $c:{\rm Gr}({A})\to \overline{\mathbb{R}}$ is Borel-measurable and bounded from below;
\item[(d)] the transition probability $q$ is regular, that is, $q(\,\cdot\,|x,a)$ is a probability measure on $(\X,\mathcal{B}(\X))$ for each $(x,a)\in {\rm Gr}(A),$ and the function $q(B|x,a)$ is Borel-measurable in $(x,a)\in {\rm Gr}(A)$ for each $B\in\mathcal{B}(\X).$
\end{itemize}

The decision process proceeds as follows: at each time epoch $t = 0,1,\ldots,$ the
current state of the system, $x_t,$ is observed. A decision-maker chooses an action $a_t\in A(x_t),$ the
cost $c(x_t,a_t)$ is accrued, and the system moves to the next state $x_{t+1}$ according to $q(\,\cdot\,|x_t, a_t).$

For $t=0,1,\ldots$ let $\mathbb{H}_t=(\X\times \mathbb{A})^t\times \X$ be the
set of histories up to epoch $t$ and $ {\mathcal
B}(\mathbb{H}_t)=({\mathcal B}(\X)\otimes {\mathcal
B}(\mathbb{A}))^t\otimes {\mathcal B}(\X)$. A randomized
decision rule at epoch $t$ is a regular transition
probability $\pi_t:H_t\to \mathbb{A}$ concentrated on $A(x_t).$ 
A policy is a sequence $\pi=\{\pi_t\}_{t= 0,1,\ldots}$ of decision rules.
Moreover, $\pi$ is called nonrandomized, if each
probability measure $\pi_t(\,\cdot\,|h_t)$ is concentrated at one
point. A nonrandomized policy is called Markov, if all of
the decisions depend on the current state and time only. A Markov policy is called
stationary, if all the decisions depend on the current
state only. 
Note that a stationary policy is a measurable selector $\varphi:\X\to\A$ for $A.$ 

The Ionescu Tulcea theorem (\cite[pp.~140-141]{BS} or \cite[p. 178]{HLL}) implies that an initial state $x_0=x$ and a
policy $\pi$ define a unique probability $P_x^{\pi}$ on the set of
all trajectories $\mathbb{H}_{\infty}=(\X\times
\mathbb{A})^{\infty}$ endowed with the product of $\sigma$-field
defined by Borel $\sigma$-field of $\X$ and $\mathbb{A}.$ Let
$\mathbb{E}_x^{\pi}$ be an expectation with respect to
$P_x^{\pi}$. Let $\alpha\in [0,1]$ and $v_{0,\alpha}^{\pi}(x)=0.$  For a finite horizon $T=1,2,\ldots,$ the expected total discounted costs is defined as $v_{T,\alpha}^{\pi}(x):=\mathbb{E}_x^{\pi}\sum\limits_{t=0}^{T-1}\alpha^tc(x_t,a_t),$ $x\in\X.$ We usually write  $v_\alpha^\pi(x)$ instead of $v_{\infty,\alpha}^{\pi}(x).$
If 
$\alpha\in [0,1)$, then $v_\alpha^{\pi}(x)$ is an infinite-horizon
expected total discounted cost.  For $T=+\infty$ and
$\alpha=1$  we assume that the cost function $c$ takes nonnegative values.  Then $v_1^{\pi}(x)$ is an infinite-horizon
expected total undiscounted cost.
The average cost per unit time is defined as $
w^{\pi}(x):=\slim\limits_{T\to +\infty}\frac{1}{T}
v_{T,1}^\pi(x),$ $x\in\X.$
For any function $g^{\pi}(x)$, including
$g^{\pi}(x)=v_{T,\alpha}^{\pi}(x)$,
$g^{\pi}(x)=v_{\alpha}^{\pi}(x)$, and $g^{\pi}(x)=w^{\pi}(x)$,
define the optimal cost $g(x):=\inf\limits_{\pi\in \Pi}g^{\pi}(x),$ $x\in\X,$
where $\Pi$ is the set of all policies. A policy $\pi$ is called optimal for the respective
criterion, if $g^{\pi}(x)=g(x)$ for all $x\in \X.$

The main result of this paper for the expected total costs, Theorem~\ref{th:dcoe}, covers expected total discounted and undiscounted costs. These two criteria are broadly used in applications for finite-horizon problems.  Expected total discounted costs and average-costs per unit time are broadly used for infinite-horizon problems.  Classic applications include inventory control and control of queueing systems.  For these two classes of applications and many other problems, expected total undiscounted cots are typically infinite for all policies, and therefore the objective criterion $v_1^\pi(x)$ is not natural.  However, there are several important applications, including optimal stopping and search problems, in which the value function $v_1(x)$ takes finite values, and the criterion of total expected undiscounted costs is natural for such problems. MDPs with expected total undiscounted costs have been  studied in the literature for positive and negative costs since  pioneering fundamental contributions \cite{Bl,St}.  The main result of this paper on expected total rewards includes the case $\alpha=1,$ and it covers negative dynamic programming, that is, MDPs with nonnegative costs.  It is well-known that discounted MDPs with bounded below costs can be reduced to undiscounted MDPs with nonnegative costs.  For countable state problems there are general results for undiscounted MDPs with expected total rewards \cite{EF1, EF2, FS}, and they imply the theory for MDPs with expected total nonnegative, nonpositive, and discounted one-step costs.  For uncountable problems such theory is not available at the present time, and some particular results can be found in \cite{EF92,SS}.

The following assumption is used in this paper to prove the existence of optimal policies
 \vskip
0.9 ex

\noindent\textbf{Assumption ${\rm S}^*$}.
\begin{itemize}
\item[(i)] the function $a\mapsto c(x,a)$ is inf-compact on $A(x)$ for each $x\in\X;$
\item[(ii)] for each $x\in \X$ the transition probability  $q(\,\cdot\,|x,a)$ is setwise
continuous in  $a\in A(x),$ that is, for every bounded
measurable function $f: \X\to \mathbb{R},$ the function $\int_\X  f(z)q(dz|x,a)$ is  continuous in $a\in A(x)$ for each $x\in\X.$
\end{itemize}

Corollary~\ref{corsel:2} provides sufficient conditions  for the existence of a measurable selector stated in assumption~(b).  For example, according to Corollary~\ref{corsel:2}(iii),    such selector exists under Assumption~${\rm S}^*$(i),   if for each $x\in\X$ there exists $a\in A(x)$  with $c(x,a)<+\infty.$  In view of Theorem~\ref{th2}, the set of $x\in\X,$ for which $c(x,a)<+\infty$ for some $x\in\X,$ is measurable.  Therefore, from a modeling point of view, it is possible to merge all the states $x\in\X,$ for which $c(x,a)=+\infty$ for all $a\in A(x),$ into a single state $x^*,$ for which the action set $A(x^*)$ is a singleton $\{a^*\}$  such that $c(x^*,a^*)=+\infty$ and $p(\{x^*\}|x^*,a^*)=1.$ In view of Corollary~\ref{corsel:2}(ii), Assumption~${\rm S}^*$(i) holds for this MDP if it holds for the original MDP.

Without loss of generality, we assume that the function $c\ge 0.$  This is true for finite-horizon problems, for infinite-horizon discounted costs with discount factors less than 1, and for infinite-horizon problems with average costs per unit time; see \cite{FKZ2012} for details.  For infinite-horizon problems with expected total undiscounted costs, it is already assumed in the definition of the objective function that $c\ge 0$. For a Borel space $\S,$ let $\M(\S)$ be the set of all Borel nonnegative measurable functions $f:\S\to\overline{\mathbb{R}}.$ For any $\alpha\ge 0$ and $w\in \M(\X),$ we define
\begin{equation}\label{e:defeta}
\eta_w^\alpha(x,a)=c(x,a)+\alpha\int_\X  w(z)q(dz| x,a),\qquad\qquad (x,a)\in {\rm Gr}(A).
\end{equation}

\subsection{Expected Total   Costs}\label{subs:1}
The following theorem states  basic properties of MDPs with the expected total costs: the validity of optimality equations, existence of stationary and Markov
optimal policies, description of  sets of stationary and Markov optimal
policies, and convergence of value iterations. \cite[Theorem 2]{FKZ2012} {provides similar results for discounted MDPs} with weakly continuous transition probabilities.  {There are also interesting results on convergence of value iterations for MDPs with possibly discontinuous transition probabilities and one-step costs; see \cite{Yu} and references therein.}

\begin{theorem}\label{th:dcoe} Let Assumption~${\rm S^*}$ hold.
Then
\begin{itemize}
\item[{\rm(i)}] the functions $v_\alpha(x)$ and $ v_{t,\alpha}(x), $ $t=0,1,\ldots,$ belong to the set $\M(\X\times[0,1]),$ and $v_{t,\alpha}(x)\uparrow
v_\alpha (x)$ as $t \to \infty$ for all $(x,\alpha)\in \X\times[0,1];$
\item[{\rm(ii)}] for each $x\in\X$ the functions $\alpha\mapsto v_\alpha(x)$ and $\alpha\mapsto v_{t,\alpha}(x),$ $t={1},{2},\ldots,$ where $\alpha\in [0,1],$ are nondecreasing {and} {l.s.c.};
\item[{\rm(iii)}] if $t=0,1,\ldots,$ $\alpha\in[0,1],$ and $x\in\X,$ then $v_{t+1,\alpha}(x)=\min\limits_{a\in A(x)}\eta_{v_{t,\alpha}}^\alpha(x,a),$ and the nonempty sets
$A_{t,\alpha}(x):=\{a\in
A(x):\,v_{t+1,\alpha}(x)=\eta_{v_{t,\alpha}}^\alpha(x,a) \}$ satisfy the properties: {\rm(a)}~${\{(x,\alpha,a)\,:\,a\in A_{t,\alpha}(x)\}}\in\mathcal{B}(\X\times [0,1]\times\A)$, and {\rm(b)}
$A_{t,\alpha}(x)=A(x),$ if $v_{t+1,\alpha}(x)=+\infty$,  and  $A_{t,\alpha}(x)$ is compact if
$v_{t+1,\alpha}(x)<+\infty;$
\item[{\rm(iv)}] for  $T=1,2,\ldots$ and $\alpha\in[0,1]$, if for a $T$-horizon Markov policy $(\phi_0,\ldots,\phi_{T-1})$ the inclusions $\phi_{T-1-t}(x)\in A_{t,\alpha}(x)$  hold for all $x\in\X$ and for all
$t=0,\ldots,T-1,$  then this policy is $T$-horizon optimal for the discount factor $\alpha$, and, in addition,  there exist Markov optimal
$T$-horizon policies $(\phi_0^\alpha,\ldots,\phi_{T-1}^\alpha)$ for the discount factor $\alpha$ such $\phi_t^\alpha(x):\X\times[0,1]\to \A$  is  Borel measurable  for each $t=0,\ldots,T-1;$
\item[{\rm(v)}] if $\alpha\in [0,1]$ and $x\in\X,$ then $v_{\alpha}(x)=\min\limits_{a\in A(x)}\eta_{v_{\alpha}}^\alpha(x,a),$
and the nonempty sets $A_{\alpha}(x):=\{a\in
A(x):\,v_{\alpha}(x)=\eta_{v_{\alpha}}^\alpha(x,a) \}$
satisfy the properties: {\rm(a)} ${\{(x,\alpha,a)\,:\,a\in A_{\alpha}(x)\}}\in\mathcal{B}(\X\times[0,1]\times \A),$ and {\rm(b)}
 $A_{\alpha}(x)=A(x), $ if $v_{\alpha}(x)=+\infty,$ and
$A_{\alpha}(x)$ is compact if $v_{\alpha}(x)<+\infty;$
\item[{\rm(vi)}] for a discount factor  $\alpha \in [0,1],$ a stationary policy $\phi$ is optimal for an infinite-horizon problem with this discount factor if and only if $\phi(x)\in A_\alpha(x)$ for all
$x\in \X,$ and for each  $\alpha \in [0,1]$ there exists  a stationary policy $\phi^\alpha,$ which is optimal for an infinite-horizon problem with the  discount factor $\alpha,$ and $\phi_\alpha(x):\X\times[0,1]\to\A$ is a Borel measurable mapping.
\end{itemize}
\end{theorem}

Before the proof of Theorem~\ref{th:dcoe}, we provide Lemma~\ref{lm1}, which is useful for establishing continuity properties of the value functions $v_{t,\alpha}$ and
$v_\alpha(x).$ The proof of this lemma uses Theorem~\ref{th2}.
For each $(x,\alpha)\mapsto w_\alpha(x)$ from $\M(\X\times\R_+)$
we consider the function $(x,\alpha)\mapsto w_\alpha^*(x):=\inf\limits_{a\in A(x)}\eta_{w_\alpha}^\alpha(x,a)$ on $\X\times\R_+.$ We recall that, according to \cite[Definition 1.1]{FKZ1}, for $x\in\X$ a function $(\alpha,a)\mapsto f_\alpha(x,a),$ mapping $\R_+\times A(x)$ to $\RR,$ is $\K$-inf-compact on $\R_+\times A(x),$ if for each compact set $K\subset\R_+$ this function is inf-compact on $K\times A(x).$


\begin{lemma}\label{lm1}
Let Assumption~${\rm S^*}$ holds, and let $(x,\alpha)\mapsto w_\alpha(x)$ be a function from $\M(\X\times\R_+)$ such that for each $x\in\X$ the function $\alpha\mapsto w_\alpha(x)$ is nondecreasing {and} {l.s.c.} Then:
\begin{itemize}
\item[{\rm(i)}]
the function $(x,a,\alpha)\mapsto \eta_{w_\alpha}^\alpha(x,a)$  belongs to $\M({\rm Gr}(A)\times\R_+),$  it is inf-compact in $a$ on $A(x)$ for each $x\in\X$ and $\alpha\ge0,$ and for each $x\in\X$ the function $(\alpha,a)\mapsto \eta_{w_\alpha}^\alpha(x,a)$ is $\K$-inf-compact on $\R_+\times A(x);$
\item[{\rm(ii)}] for each $(x,a)\in{\rm Gr}(A)$ the function $\alpha\mapsto \eta_{w_\alpha}^\alpha(x,a)$ is nondecreasing;
\item[{\rm(iii)}]
the function $(x,\alpha)\mapsto w_\alpha^*(x)$  belongs to $\M(\X\times\R_+);$
\item[{\rm(iv)}] for each $x\in\X$ the function $\alpha\mapsto w_\alpha^*(x)$ is nondecreasing {and l.s.c.};
\item[{\rm(v)}] there exists a Borel mapping $(x,\alpha)\mapsto f_\alpha(x)$ {from} $\X\times\R_+$ into $\A$ such that $f_\alpha(x)\in A(x)$ and $w_\alpha^*(x)=\eta_{w_\alpha}^\alpha(x,f_\alpha(x))$ for all  $x\in \X$ and $\alpha\ge0;$
\item[{\rm(vi)}] the nonempty sets $A^*_\alpha(x)=\left\{a\in A(x):\,w_\alpha^*(x)=\eta_{w_\alpha}^\alpha(x,a)\right\},$ $(x,\alpha)\in \X\times\R_+,$ satisfy the following properties: {\rm(a)}~${\{(x,\alpha,a)\,:\,a\in A_\alpha^*(x)\}}\in\mathcal{B}(\X\times \R_+\times\mathbb{A});$ {\rm(b)}
 $A^*_\alpha(x)=A(x)$, if $w_\alpha^*(x)=+\infty,$ and $A_\alpha^*(x)$ is compact  if $w_\alpha^*(x)<+\infty.$
\end{itemize}
\end{lemma}
\begin{proof}
Let us prove (i{,ii}).
The function $(x,a,\alpha)\mapsto \eta_{w_\alpha}^\alpha(x,a)$ is nonnegative and nondecreasing in $\alpha$ because $(x,a)\mapsto c(x,a)$ and $(x,\alpha)\mapsto w_\alpha(x)$ are nonnegative and nondecreasing in $\alpha.$ Borel-measurability and continuity properties of $(x,\alpha)\mapsto w_\alpha(x)$ and Assumptions (a) and (d) imply that the function $(x,a,\alpha)\mapsto \int_\X  w_\alpha(z)q(dz|x,a)$ is Borel  measurable on ${\rm Gr}(A)\times\R_+,$ which implies that the function $(x,a,\alpha)\mapsto \eta_{w_\alpha}^\alpha(x,a)$ is Borel measurable on ${\rm Gr}(A)\times\R_+.$ The function $a\mapsto \alpha\int_\X  w_\alpha(z)q(dz|x,a)$ is l.s.c. on $A(x)$ for each $x\in\X$ and $\alpha\ge0.$ This follows from Assumption~${\rm S^*},$ and \cite[Theorem 4.2]{FKL}. For each $(x,\alpha)\in\X\times\R_+$ the function $a\mapsto \eta_{w_\alpha}^\alpha(x,a)$ is inf-compact on $A(x)$ as the sum of an inf-compact function and a nonnegative l.s.c. functions. Let us fix an arbitrary $x\in\X$. The function $c(x,\cdot)$ is inf-compact in view of Assumption~${\rm S^*}$(i) and does not depend on $\alpha.$  Therefore, the function $(\alpha,a)\mapsto c(x,a)$ is $\K$-inf-compact on $\R_+\times A(x).$  As follows from \cite[Theorem 4.2]{FKL}, the function $(\alpha,a)\mapsto \alpha\int_\X  w_\alpha(z)q(dz|x,a)$ is l.s.c. on $\R_+\times A(x).$ Moreover, this function is nonnegative. The function $(\alpha,a)\mapsto \eta_{w_\alpha}^\alpha(x,a)$ is $\K$-inf-compact on $\R_+\times A(x)$ because it is a sum of a $\K$-inf-compact function and a nonnegative l.s.c. function.

Let us prove statements~(iii,v,vi). Statement~(i), Theorem~\ref{th2}, and Corollary~\ref{corsel:2}(i) directly imply statements~(iii) and (v). Property (vi)(a) follows from Borel measurability of $(x,a,\alpha)\mapsto\eta_{w_\alpha}^\alpha(x,a)$ on ${\rm Gr}(A)\times\R_+$ and $(x,\alpha)\mapsto w^*_\alpha(x)$ on $\X\times\R_+;$ and property (vi)(b) follows from inf-compactness of $a\mapsto \eta_{w_\alpha}^\alpha(x,a)$ on $A(x)$ for each $(x,\alpha)\in\X\times\R_+.$

Let us prove statement~(iv). According to statement (i), the function $(\alpha,a)\mapsto \eta_{w_\alpha}^\alpha(x,a)$ is $\K$-inf-compact on $\R_+\times A(x).$
Therefore, in view of Berge's theorem for noncompact action sets \cite[Theorem 1.2]{FKZ1}, the function $\alpha\mapsto w_\alpha^*(x)$ is l.s.c. on $\R_+.$
\end{proof}

\begin{proof}[Proof of Theorem~\ref{th:dcoe}]
According to \cite[Proposition 8.2]{BS}, the functions $v_{t,\alpha}(x),$ $t=0,1,\ldots,$
recursively satisfy the optimality equations
with $v_{0,\alpha}(x)=0$ and
$v_{t+1,\alpha}(x)=\inf\limits_{a\in A(x)}\eta_{v_{t,\alpha}}^\alpha(x,a),$ for all $(x,\alpha)\in
\X\times[0,1].$ So, Lemma~{\ref{lm1}}(i) sequentially applied to the functions $v_{0,\alpha}(x),$ $v_{1,\alpha}(x),\ldots,$ implies statement{~(i)} for them,  {and Lemma~\ref{lm1}(iv)} implies that these functions are {l.s.c. in $\alpha.$} According to \cite[Proposition~9.17]{BS}, $v_{t,\alpha}(x)\uparrow
v_{\alpha}(x)$ as $t \to +\infty$ for each $(x,\alpha)\in\X\times[0,1].$ Therefore, $v_\alpha(x)\in\M(\X\times[0,1]),$ and $v_\alpha(x)$ is nondecreasing and l.s.c. in $\alpha{.}$ Thus, statement{s (i,ii) are} proved.

In addition, \cite[Lemma 8.7]{BS} impl{ies} that a Markov policy defined at the first $T$ steps by the mappings $\phi_0^\alpha,...\phi_{T-1}^\alpha,$ that satisfy for all
$t=1,\ldots,T$ the equations $v_{t,\alpha}(x)=\eta_{v_{t-1,\alpha}}^\alpha(x,\phi_{T-t}^\alpha(x)),$ for each $(x,\alpha)\in \mathbb{X}\times[0,1],$ is optimal for the horizon $T.$ According to \cite[Propositions~9.8 and 9.12]{BS}, $v_{\alpha}(x)$ satisfies the discounted cost optimality equation $v_{\alpha}(x)=\inf\limits_{a\in A(x)}\eta_{v_{\alpha}}^\alpha(x,a)$ for each $(x,\alpha)\in \X\times[0,1];$ and a stationary policy $\phi_\alpha(x)$ is discount-optimal if and
only if $v_{\alpha}(x)=\eta_{v_{\alpha}}^\alpha(x,\phi_\alpha(x))$ for each $x\in \X.$  Statements  (iii-vi) follow from these facts and Lemma~\ref{lm1} (v,vi).
\end{proof}

\subsection{Average Costs Per Unit Time}\label{subs:2}

Following \cite{Schal}, we assume that $ w^*:=\inf\limits_{x\in
\X}w(x)<+\infty,$ that is, there exist $x\in\X$ and
$\pi\in\Pi$ with $w^\pi(x)<+\infty.$ Otherwise, if this assumption
does not hold, then the problem is trivial, because $w(x)=+\infty$
for all $x\in\X$ and any policy $\pi$ is average-cost optimal.

Define the following quantities for $\alpha\in [0,1)$:
\begin{equation*}
m_{\alpha}=\inf\limits_{x\in \X}v_{\alpha}(x),\quad
u_{\alpha}(x)=v_{\alpha}(x)-m_{\alpha},
\end{equation*}
\begin{equation*}
\underline{w}=\ilim\limits_{\alpha\uparrow
1}(1-\alpha)m_{\alpha},\quad\overline{w}=\slim\limits_{\alpha\uparrow
1}(1-\alpha)m_{\alpha}.
\end{equation*}
According to \cite[Lemma 1.2]{Schal},
\begin{equation}\label{eq:schal} 0\le \underline{w}\le \overline{w}\le w^*<
+\infty.
\end{equation}


In this section we show that Assumption~${\rm S^*}$ and
boundedness assumption Assumption~${\rm \underline{B}}$ on the
function $u_\alpha$ introduced in \cite{FKZ2012}, which is weaker than  boundedness
Assumption~${\rm B}$ introduced in~\cite{Schal}, lead
to the validity of stationary average-cost optimal inequalities and
the existence of stationary policies. Stronger results hold under stronger
Assumption~${\rm B};$ see \cite{HL} and the conclusions of \cite[Theorem 4]{FKZ2012}. Example~\ref{exa} in Section~\ref{s:exa} describes an MDP satisfying Assumption~${\rm \underline{B}}$ and not satisfying Assumption~${\rm B}$.

\textbf{Assumption~${\rm \underline{B}}$.} $u(x):=\ilim\limits_{\alpha \uparrow
1}u_{\alpha}(x)<+\infty$ for all $x\in \X$.

\textbf{Assumption~${\rm B}$.} $\sup_{\alpha\in [0,1)}u_{\alpha}(x)<+\infty$
for all $x\in \X$.
\vskip
0.9 ex

In the rest of this paper we assume that Assumption~${\rm \underline{B}}$ holds.
In view of Theorem~\ref{th:dcoe}(i), if $v_\alpha(x)=+\infty$ for some $(x,\alpha)\in\X\times [0,1),$ then $u_\beta(x)=v_\beta(x)=+\infty$ for all $\beta\in [\alpha,1),$ and $u(x)=+\infty,$ where $m_{\beta}$ is finite in view of \eqref{eq:schal}.  Thus  Assumption~${\rm \underline{B}}$ implies that  $v_\alpha(x)<+\infty,$  and therefore $u_\alpha(x)<+\infty$ for all $(x,\alpha)\in\X\times [0,1).$ {Then the function $\alpha\mapsto m_{\alpha+}:=\lim_{\beta\downarrow \alpha} m_\beta$ is real-valued on $[0,1).$ Moreover, it is nondecreasing on $[0,1)$ because $c\ge0.$ Thus, this function is u.s.c. on $[0,1)$ by its definition}.  Therefore,
${\tilde{u}}_\alpha(x)=v_\alpha(x)-m_{\alpha{+}}$ {is Borel measurable on $\X\times [0,1),$} and this function is l.s.c. in $\alpha\in [0,1)$ for each $x\in \X.$

Let us define the following nonnegative functions on $\X,$
\begin{equation}\label{eq:defuuU}
U_\beta(x) :=  \inf\limits_{\alpha\in \left[\beta,1
\right)}u_{\alpha}(x), \quad \beta\in [0,1),\
x\in \X.
\end{equation}
Under Assumption~${\rm S^*},$ for each $\beta\in [0,1)$ the function $U_\beta:\X\to\R_+$ is Borel measurable. {Indeed, let us consider a sequence $\beta_n\uparrow 1$ with $\beta_1:=\beta.$  Then $U_\beta(x)=\inf_{n\ge 1}  \inf\limits_{\alpha\in \left[\beta_n,\beta_{n+1}
\right)}u_{\alpha}(x).$ Therefore, the Borel measurability of the functions $U_\beta(\cdot)$ follows from the Borel measurability on $\X$ of the functions
\[
U_{\beta,\bar{\beta}}(x) :=  \inf\limits_{\alpha\in \left[\beta,\bar{\beta}
\right)}u_{\alpha}(x), \quad
x\in \X.
\]
defined for all $\beta\in [0,1)$ and for all $\bar{\beta}\in (\beta,1).$ We follow the convention $\inf\emptyset:=+\infty. $

Let us prove that the function $U_{\beta,\bar{\beta}}(x)$ is Borel measurable on $\X$ for all fixed $\beta\in [0,1)$ and $\bar{\beta}\in (\beta,1).$ The function $m_\alpha$ is nondecreasing and bounded on $[\beta,\bar{\beta}).$  Thus, it can have only positive jumps whose total sum is finite. The set $D:=\{\alpha\in [\beta,\bar{\beta}): m_\alpha\ne m_{\alpha+}\}$ is countable or finite. Moreover, for each $\varepsilon>0$ the number of jumps which are larger than $\varepsilon$ is finite. Therefore, let us enumerate its elements, $D=\{\alpha_1,\alpha_2,\ldots\},$ where $D$ is finite or countable, in a way that the jumps $(m_{\alpha_n+}- m_{\alpha_n-})$ do not increase in  $n,$ that is, $(m_{\alpha_n+}- m_{\alpha_n-})\ge (m_{\alpha_{n+1}+}- m_{\alpha_{n+1}-}),$ if $D$ consists of more than $n$ points, where $ m_{\alpha-}:=\lim_{\tilde{\alpha}\uparrow \alpha} m_{\tilde{\alpha}}.$ Then $\lim_{n\to\infty} (m_{\alpha_n+}- m_{\alpha_n-})=0$ if the set $D$ is infinite. So, in this case, for a fixed $\varepsilon>0$ there exists $n(\varepsilon)=1,2,\ldots$ such that $(m_{\alpha_k+}- m_{\alpha_k-})\le \varepsilon$ for each $k=n(\varepsilon)+1,n(\varepsilon)+2,\ldots,$ and we set $D_{n(\varepsilon)}:=\{\alpha_1,\alpha_2,\ldots,\alpha_{n(\varepsilon)}\}.$ When $D$ is finite, then $D = D_{n(\varepsilon)}$ for each $\varepsilon \in(0,\frac1k),$ where $k=1,2,\ldots$ is sufficiently large. Theorem~\ref{th1} applied to the Borel spaces $\X$ and $[\beta,\bar{\beta}),$ the set-valued map $B(x)=[\beta,\bar{\beta}){\setminus D_{n(\varepsilon)}}$ for all $x\in\X,$ the function $u(x,\alpha):={\tilde{u}}_\alpha(x){:\X\times [\beta,\bar{\beta})\to \R}{,}$ and the $\sigma$-compact set $B^f(x)=\{\tilde{\beta}\in [\beta,\bar{\beta})\setminus D_{n(\varepsilon)}\,:\,\tilde{u}_{\tilde{\beta}}(x) <+\infty\}=[\beta,\bar{\beta})\setminus D_{n(\varepsilon)}$ implies that the function $x\mapsto\inf\limits_{\alpha\in \left[\beta,\bar{\beta}
\right)\setminus D_{n(\varepsilon)}}\tilde{u}_{\alpha}(x)$ is Borel measurable for each $\varepsilon>0.$ Therefore, the function $x\mapsto\inf\limits_{\alpha\in \left[\beta,\bar{\beta}
\right)\setminus D}\tilde{u}_{\alpha}(x)$ is Borel measurable as a pointwise monotone limit of the sequence of Borel measurable functions $\{x\mapsto\inf\limits_{\alpha\in \left[\beta,\bar{\beta}
\right)\setminus D_{n(\frac1k)}}\tilde{u}_{\alpha}(x)\,:\,k=1,2,\ldots\}.$ This convergence follows from
\[\inf\limits_{\alpha\in \left[\beta,\bar{\beta}
\right)\setminus D_{n(\frac1k)}}\tilde{u}_{\alpha}(x)\le \inf\limits_{\alpha\in \left[\beta,\bar{\beta}
\right)\setminus D}\tilde{u}_{\alpha}(x)\le \inf\limits_{\alpha\in \left[\beta,\bar{\beta}
\right)\setminus D_{n(\frac1k)}}\tilde{u}_{\alpha}(x)+\frac{1}{k},\quad k=1,2,\ldots,
\]
where the second inequality follows from the left-continuity in $\alpha$ of the function $v_\alpha(x)-m_{\alpha-}$ on $(0,1).$  The function $x\mapsto U_{\beta,\bar{\beta}}(x)=\min\{\inf\limits_{\alpha\in \left[\beta,\bar{\beta}
\right)\cap D}u_{\alpha}(x),\inf\limits_{\alpha\in \left[\beta,\bar{\beta}
\right)\setminus D}\tilde{u}_{\alpha}(x)\}$ is Borel measurable on $\X$ as a minimum of two Borel measurable functions, where the function $x\mapsto\inf\limits_{\alpha\in \left[\beta,\bar{\beta}
\right)\cap D}u_{\alpha}(x)$ is Borel measurable on $\X$ as an infumum of the at most countable family of functions $\{x\mapsto u_{\alpha}(x)\,:\,\alpha\in \left[\beta,\bar{\beta}
\right)\cap D\}$ which are Borel measurable on $\X$ according to Theorem~\ref{th:dcoe}(i).}


In view of the definition of $u$ in Assumption~${\rm \underline{B}},$
\begin{equation}\label{eq5.5} {
u(x)=\lim\limits_{\alpha\uparrow
1}U_{\alpha}(x),} \qquad x\in\X,
\end{equation}
{this equality implies that} the function $u$ is Borel measurable, if Assumptions~${\rm S^*}$ holds, and under this assumption the following sets can be defined for
$u$ introduced in Assumption ${\rm \underline{B}}$:
\[
\begin{aligned}
&A^{u}(x):=\left\{a\in A(x)\, : \,
\overline{w}+u(x)\ge \eta_u^1(x,a) \right\},\\
&A_{u}(x):=\left\{a\in A(x)\,:\,\min_{a^*\in A(x)}\eta_u^1(x,a^*)=\eta_u^1(x,a)\right\}, \quad x\in\X.
\end{aligned}
\]
In view of Lemma~\ref{lm1}, the sets $A_{u}(x)$ are nonempty and compact for all $x\in \X.$
In the following theorem
we show that Assumption~${\rm S^*}$ and
boundedness assumption Assumption ${\rm \underline{B}}$ on the
functions $\{u_\alpha\}_{\alpha\in(0,1)},$ which is weaker than the boundedness
Assumption~${\rm B}$ introduced in \cite{Schal}, lead
to the validity of stationary average-cost optimal inequalities and
the existence of stationary policies. Stronger facts under
Assumption~${\rm B}$ are established in \cite{HL}. \cite[Theorems 3 and 4]{FKZ2012} are respectively counterparts to Theorem 3.3 and the main result in \cite{HL} for MDPs with weakly continuous transition probabilities.  Assumption~${\rm B}$ and some additional conditions lead to the validity of optimality equations for average-costs MDPs. In \cite{FLi17} such sufficient conditions are provided for MDPs with weakly continuous transition probabilities and applied to inventory control. More general sufficient conditions for validity of optimality equations are provided in \cite[Section 7]{FKL} for MDPs with weakly and setwise continuous transition probabilities.

\begin{theorem}\label{teor1}
Let Assumptions~${\rm S^*}$ and ${\rm \underline{B}}$ hold. Then for infinite-horizon average costs per unit time there exists a stationary optimal policy $\phi$ satisfying
\begin{equation}\label{eq7111}
\overline{w}+u(x)\ge \eta_u^1(x, \phi(x)),\quad x\in \X,
\end{equation}
with $u$ defined in Assumption ${\rm \underline{B}}$, and for this policy 
\begin{equation}\label{eq:7121}
w(x)=w^{\phi}(x)=\slim\limits_{\alpha\uparrow
1}(1-\alpha)v_{\alpha}(x)=\overline{w}=w^*,\quad x\in \X.
\end{equation}
 Moreover, the following statements hold:
\begin{itemize}
\item[{\rm(a)}] the function $u:\X\to \R_+$ defined {in} Assumption ${\rm \underline{B}}$ is Borel measurable;
\item[{\rm(b)}] the nonempty sets $A^{u}(x)$, $x\in\X$,
satisfy the following properties: ${\rm(b_1)}$ ${\rm Gr}(A^{u})\in\mathcal{B}(\X\times \mathbb{A});$ ${\rm(b_2)}$ for each $x\in\X$ the set $A^{u}(x)$ is compact;
\item[{\rm(c)}] if $\varphi(x)\in A^{u}(x)$ for all $x\in \X$ for a stationary policy $\varphi,$ then $\varphi$   satisfies \eqref{eq7111} and \eqref{eq:7121}, with $u$ defined in
Assumption ${\rm \underline{B}}$ and with $\phi=\varphi,$ and $\varphi$ is optimal for
average costs per unit time;
\item[{\rm(d)}] the sets $A_u(x)$ are compact and $A_u(x)\subset A^u(x)$ for all $x\in X,$ and there exists a stationary policy $\varphi$ with
$\varphi(x)\in A_{u}(x)\subset
A^{u}(x)$ for all $x\in\X.$
\end{itemize}
\end{theorem}
The proof of Theorem~\ref{teor1} uses the following statement.
\begin{lemma}\label{lemma3} Under Assumptions~${\rm \underline{B}}$ and ${\rm S^*}$, 
\begin{equation}\label{eq:l31} \overline{w}+u(x)\ge
\min\limits_{a\in A(x)} \eta_u^1(x,a)
,\quad x\in\X.
\end{equation}
\end{lemma}
\begin{proof}
Fix an arbitrary $\varepsilon^*>0$. Due to the definition of $\overline{w}$, there
exists $\alpha_0\in (0,1)$ such that
\begin{equation}\label{eq:l32}
\overline{w} +\varepsilon^* > (1-\alpha)m_{\alpha},
\quad \alpha\in [\alpha_0,1).
\end{equation}
The $\R_+$-valued function $U_\alpha$ is Borel measurable for all ${\alpha\in(0,1)}.$ Therefore, the
function $\eta_{U_\alpha}^\alpha(x,a)$ is well-defined.
Let us prove that
\begin{equation}\label{eq:l33}
\overline{w}+\varepsilon^*+u(x)\ge\min\limits_{a\in A(x)}
\eta_{U_\alpha}^\alpha(x,a),\quad x\in\X,\, \alpha\in [\alpha_0,1).
\end{equation}

Indeed, Theorem~\ref{th:dcoe}(v) and (\ref{eq:l32}) imply that
\[
\begin{aligned}
&\overline{w}+\varepsilon^*+u_{\beta}(x)>
(1-\beta)m_\beta+u_{\beta}(x)=v_\beta(x)-\beta m_\beta
\\ &=\min\limits_{a\in
A(x)}\eta_{u_\beta}^\beta(x,a)
\ge \min\limits_{a\in A(x)}\eta_{U_\alpha}^\alpha(x,a),
\end{aligned}
\]
for each $x\in\X$ and
$\alpha,\beta\in [\alpha_0,1)$ such that $\beta\ge \alpha.$ Since the right-hand side of the above inequality does not depend on $\beta\in[\alpha,1),$ by taking the infimum in $\beta\in [\alpha,1),$ we
obtain 
\[
\overline{w}+\varepsilon^*+U_\alpha(x)\ge \min\limits_{a\in A(x)}\eta_{U_\alpha}^\alpha(x,a),
\]
for all $x\in\X$ and $\alpha\in [\alpha_0,1).$ Therefore, since the function $U_\alpha(x)$ is nonincreasing in $\alpha\in (0,1),$ inequalities \eqref{eq:l33} hold in view of \eqref{eq5.5}.

Let us fix an arbitrary $x\in\X$ {and $\alpha\in [\alpha_0,1).$} By Lemma~\ref{lm1}{(vi) applied to $w_\beta(y)= U_\alpha(y),$ $y\in\X,$ $\beta \in[0,1),$} there exists
$a_\alpha\in A(x)$ such that $ \min\limits_{a\in
A(x)}\eta_{U_\alpha}^{\alpha}(x,a)=
\eta_{U_\alpha}^{\alpha}(x,a_{\alpha}).$
 Since $U_\alpha\ge 0$, for
$\alpha\in [\alpha_0,1),$  inequality  (\ref{eq:l33}) can be
continued as
\begin{equation}\label{eq:5.14}
\overline{w}+\varepsilon^*+u(x)\ge
\eta_{U_\alpha}^{\alpha}(x,a_{\alpha})\ge
c(x,a_{\alpha}).
\end{equation}
Thus, for all $\alpha\in [\alpha_0,1)$ {and for level sets $\mathcal{D}$ defined in the beginning of Section~\ref{s2},}
\[
a_\alpha\in
\mathcal{D}_{\eta_{U_\alpha}^{\alpha}(x,\,\cdot\,)}(\overline{w}+\varepsilon^*+u(x))\subset
\mathcal{D}_{c(x,\,\cdot\,)}(\overline{w}+\varepsilon^*+u(x))\subset
A(x).
\]
Since the function $c(x,\,\cdot\,)$ is inf-compact, the nonempty set
$\mathcal{D}_{c(x,\,\cdot\,)}(\overline{w}+\varepsilon^*+u(x))$ is
compact. Therefore, for every sequence $\beta_n\uparrow 1$ of numbers
from $[\alpha_0,1)$ there is a  subsequence $\{\alpha_n\}_{n\ge
1}$ such that the sequence $\{a_{\alpha_n}\}_{n\ge 1}$ converges
and $a_*:=\lim_{n\to\infty} a_{\alpha_n}\in A(x)$.
Consider a sequence $\alpha_n\uparrow 1$ such that
$a_{\alpha_n}\to a_*$ for some $a_*\in A(x).$
 Due to \cite[Corollary~4.2]{FKL} and \eqref{eq5.5},
\[
\ilim\limits_{n\to \infty}
\alpha_{n}\int_\X
U_{\alpha_n}(z)q(dz|x,a_{n})\ge
\int_\X u(z)q(dz|x,a_*).
\]
Therefore, since the function $c$ is lower semi-continuous in $a,$ we have that (\ref{eq:5.14})
 impl{ies}
\[
\begin{aligned}
&\overline{w}+\varepsilon^*+u(x)\ge \limsup\limits_{n\to\infty}
\eta_{U_{\alpha_n}}^{\alpha_n}(x,a_{\alpha_n})
\\
&\ge
c(x,a_*)+\int_\X u(z)q(dz|x,a_*)\ge\min_{a\in A(x)} \eta_u^1(x,{a}),
\end{aligned}
\]
which implies \eqref{eq:l31} because $\varepsilon^*>0$ is  arbitrary.
\end{proof}

\begin{proof}[Proof of Theorem~\ref{teor1}]For statement
(a) see \eqref{eq5.5} and the following sentence. Since ${\rm
Gr}(A^{u})=\{(x,a)\in {\rm Gr}(A):\, g(x,a)\ge 0\}$, where
$g(x,a)=\overline{w}+u(x)- c(x,a)-\int_\X  u(y)q(dy|x,a)$ is a Borel
function, the set ${\rm Gr}(A^{u}) $ is Borel. The sets $A^{u}(x)$,
$x\in\X$, are compact because for each $x\in\X$ the function $a\mapsto \eta_{u}^1(x,a)$ is inf-compact on $A(x)$ as a sum of inf-compact and nonnegative l.s.c. functions. Thus, statement (b) is proved.  The Arsenin-Kunugui theorem
implies the existence of a stationary policy $\phi$ such that
$\phi(x)\in A^{u}(x)$ for all $x\in\X.$ Statement (d)
follows from and Lemma~\ref{lm1}(v) because each $a_*\in A_u(x)$ satisfies $\eta_u^1(x,a_*)=\min_{a^*\in A(x)}\eta_u^1(x,a^*)\le\overline{w}+u(x),$ where the inequality holds since $A^u(x)\ne \emptyset.$ The
remaining conclusions of Theorem~\ref{teor1} follow from Lemma~\ref{lemma3} and \cite[Theorem~1]{FKZ2012} stating that inequalities \eqref{eq7111}
imply optimality of the policy $\phi$ and \eqref{eq:7121}.
\end{proof}

\section{An Example Showing that Assumption~${\rm B}$ is Stronger than Assumption~${\rm \underline{B}}$ }\label{s:exa}
This section presents an example of an MDP satisfying Assumption~${\rm \underline{B}}$ and not satisfying Assumption~${\rm B}.$  This MDP has a countable state space $\X$ and a  decision space $\A$ consisting of a single point. Since ${\A}$ is a singleton, this MDP is defined by a countable state set of states $\X,$ transition probabilities $q(y|x),$ and one-step costs $c(x)\ge 0,$ where $x,y\in\X.$ Any such MDP satisfies Assumption ${{\rm S^*}.}$  If the discrete topology is considered on $\X,$ any such MDP also satisfies the analogous continuity condition ${W^*}$ for MDPs with weakly continuous transition probabilities studied in \cite{FKZ2012}.  The question whether Assumption~${\rm \underline{B}}$ is indeed stronger than Assumption~${\rm B}$ remained open since Assumption~${\rm \underline{B}}$ was introduced in \cite{FKZ2012}.  An earlier attempt to answer this question lead to constructing in \cite{BFZ} some nontrivial sequences related to the Tauberian and Hardy-Littlewood theorems.

For two numbers $\beta\in(0,1)$ and $M>0,$ let  \[\begin{aligned}&\eb:=1-\beta>0,\\ &\gbm:=\max\{\frac{\beta+1}{2},1-\frac{\eb}{3M}\}\in(\beta,1),\\
&\nbm:=\lfloor\log_\gbm(\min\{\frac12,\frac{M(1-\gbm)}{\eb}\}){\rfloor}+1\ge1,\\ &\dbm:=\max\{\frac{\gbm+1}{2},(1-\frac{1}{\eb\nbm})^{\frac{1}{\nbm}}\}\in (\gbm,1),\end{aligned}\] where $\lfloor \cdot\rfloor$ is the integer part of a number, and for $\alpha\in(0,1)$ let \[\pbm(\alpha):=\eb\frac{(1-\alpha^\nbm)^2}{1-\alpha}.\]

Parameters of the example are generated by the following procedure:
\begin{itemize}
\item choose an arbitrary $\alpha^{(1)}\in [\frac{1}{2},1);$
\item for $n=1,2,\ldots$ set
$\varepsilon^{(n)}:=\varepsilon_{\alpha^{(n)}},$ $\gamma^{(n)}:=\gamma_{\alpha^{(n)},n},$ $N{(n)}:=n^*_{\alpha^{(n)},n},$ and $\alpha^{(n+1)}:=\delta_{\alpha^{(n)},{n}}.$
\end{itemize}
 Note that $1-\frac{1}{2^n}\le\alpha^{(n)}<\gamma^{(n)}< \alpha^{(n+1)}<1$ for each $n=1,2,\ldots,$ where the first inequality follows from $\alpha^{(n+1)}\ge \frac{1+\alpha^{(n)}}{2}$ and $\alpha^{(1)}\ge \frac12.$ Thus, $\alpha^{(n)}\uparrow 1$ and $\gamma^{(n)}\uparrow 1$ as $n\to +\infty.$
 \begin{example}\label{exa}
{\rm Consider an MDP defined by the state space $\X:=\{0\}\cup\{(n,k)\,:\,n=1,2,\ldots,\,k=1,2,\ldots,2N{(n)}\},$  by the transition probabilities $q(0|0):=1,$ $q(0|n,2N{(n)}):=1,$ and $q(n,k+1|n,k):=1,$ where $n=1,2,\ldots,$ and $k=1,2,\ldots,2{N(n)}-1,$ and  by one-step costs $c(0):=1,$ $c(n,k):=1-\varepsilon^{(n)},$ if $k=1,2,\ldots,N(n),$ and $c(n,k):=1+\varepsilon^{(n)}$ if $k=N(n)+1,N(n)+2,\ldots, 2N(n),$ $n=1,2,\ldots.$ This MDP has a single policy because there is only one possible action, whose notation we omit.
 }
\end{example}

\begin{proposition}\label{prop1}
The MDP defined in Example~\ref{exa} satisfies Assumption~${\rm \underline{B}},$ but it does not satisfy Assumption~${\rm B}.$
\end{proposition}

The proof of Proposition~\ref{prop1} uses the following lemma.

\begin{lemma}\label{lem:exa1} If $\beta\in (0,1)$ and $M>0,$ then
$\pbm(\alpha)\le1$ for each $\alpha\in(0,\beta]\cup[\dbm,1),$ and $\pbm(\gbm)\ge M.$
\end{lemma}
\begin{proof}
The inequality $\pbm(\alpha)\le1$ holds for  $\alpha\in(0,\beta]$ because $\eb\frac{(1-\alpha^\nbm)^2}{1-\alpha}\le \frac{1-\beta}{1-\alpha}\le1.$ The inequality $\pbm(\alpha)\le1$ holds for each $\alpha\in[\dbm,1)$ because $\alpha\ge (1-\frac{1}{\eb\nbm})^{\frac{1}{\nbm}}$ for each $\alpha\in[\dbm,1),$  according to the definition of $\dbm,$ and, therefore,
\[
1\ge\eb(1-\alpha^\nbm)\nbm \ge\eb(1-\alpha^\nbm)(1+\alpha+\ldots+\alpha^{\nbm-1})= \eb\frac{(1-\alpha^\nbm)^2}{1-\alpha}=\pbm(\alpha),
\]
where the second inequality holds since $\alpha\in (0,1).$

To prove  $\pbm(\gbm)\ge M,$ we observe that
\begin{equation}\label{eq:exa1}
\eb(1+\gbm+\ldots+(\gbm)^{\nbm-1})=\frac{\eb}{1-\gbm} - \frac{\eb(\gbm)^\nbm}{1-\gbm}\ge 3M-M=2M,
\end{equation}
where the inequality $\frac{\eb}{1-\gbm}\ge 3M$ follows from the definition of $\gbm,$ which implies that $\gbm\ge1-\frac{\eb}{3M},$ and the inequality $\frac{\eb(\gbm)^\nbm}{1-\gbm}\le M$
holds because $\nbm\ge \log_\gbm(\frac{M(1-\gbm)}{\eb})$  according to the definition of $\nbm,$ and hence, since $\gbm\in(0,1)$, we have  $(\gbm)^\nbm\le\frac{M(1-\gbm)}{\eb}.$
Moreover, the definition of $\nbm$ implies that $\nbm\ge \log_\gbm(\frac12),$ which is equivalent to $(\gbm)^\nbm\le \frac12$ because $\gbm\in(0,1).$ So, $1-(\gbm)^\nbm\ge \frac12.$ Thus, due to \eqref{eq:exa1},
\[
\pbm(\gbm) = \eb(1+\gbm+\ldots+(\gbm)^{\nbm-1})(1-(\gbm)^\nbm)\ge M,
\]
where the equality holds because $\frac{(1-(\gbm)^\nbm)}{1-\gbm}=1+\gbm+\ldots+(\gbm)^{\nbm-1}.$

\end{proof}
 \begin{proof}[Proof of Proposition~\ref{prop1}]

According to Lemma~\ref{lem:exa1}, for each $n=1,2,\ldots$
\begin{equation}\label{eq:exa2}
\varepsilon^{(n)}\frac{(1-\alpha^{N{(n)}})^2}{1-\alpha}\le 1\mbox{ for each } \alpha \in (0,\alpha^{(n)}]\cup[\alpha^{(n+1)},1),\mbox{ and }\varepsilon^{(n)}\frac{(1-(\gamma^{(n)})^{N{(n)}})^2}{1-\gamma^{(n)}}\ge n.
\end{equation}

Fix $\alpha\in(0,1)$ and $n=1,2,\ldots.$ Note that
\begin{equation}\label{eq:exa4}
v_\alpha(n,k)=
\left\{
\begin{array}{l}
\frac{1}{1-\alpha}+\varepsilon^{(n)}\frac{(1-\alpha^{N(n)})\alpha^{N(n)-k+1}-(1-\alpha^{N(n)-k+1})}{1-\alpha},\quad k=1,2,\ldots, N(n),\\
\frac{1}{1-\alpha}+\varepsilon^{(n)}\frac{1-\alpha^{2N(n)-k+1}}{1-\alpha},\quad k=N(n)+1,N(n)+2,\ldots , 2N(n),\\
\end{array}
\right.
\end{equation}
and $v_\alpha(0)=\frac{1}{1-\alpha}.$ Therefore, $v_\alpha(n,1)<v_\alpha(0)<v_\alpha(n,N(n)+k),$ $k=1,2,\ldots, N(n).$ In addition, $v_\alpha(n,1)\le v_\alpha(n,k)$ for $k=1,2,\ldots,N(n)$ since $v_\alpha(n,k)=\frac{1-\varepsilon^{(n)}}{1-\alpha}+\frac{\varepsilon^{(n)}\alpha^{N(n)+1}(2-\alpha^{N(n)})}{1-\alpha}\alpha^{-k}$ for these values of $k.$
Therefore, {for all $\alpha\in [0,1)$}
\begin{equation}\label{eq:exa3}
m_\alpha=\inf_{n=1,2,\ldots}v_\alpha(n,1) = \frac{1}{1-\alpha}-\sup_{n=1,2,\ldots}\varepsilon^{(n)}\frac{(1-\alpha^{N(n)})^2}{1-\alpha}. 
\end{equation}

According to \eqref{eq:exa2}--\eqref{eq:exa3}, $v_{\gamma^{(n)}}(0)-m_{\gamma^{(n)}}\ge n\to +\infty$ as $n\to+\infty.$ That is, Assumption~${\rm B}$ does not hold. On the other hand, \eqref{eq:exa2}--\eqref{eq:exa3} imply that $v_{\alpha^{(n)}}(0)-m_{\alpha^{(n)}}\le 1,$ and $v_\alpha(\tilde{n},\tilde{k})-v_\alpha(0)\le N(\tilde{n})$ for all $\alpha\in[0,1),$ $\tilde{n}=1,2,\ldots,$ $\tilde{k}=1,2,\ldots, 2N(\tilde{n}),$ and $n=1,2,\ldots.$ Thus, $\ilim_{\alpha\uparrow1}u_\alpha(x)<+\infty$ for each $x\in\X.$ Assumption~${\rm \underline{B}}$ holds.\end{proof}

%

\section{Acknowledgements} {We thank Janey (Huizhen) Yu for valuable remarks.}
The first author was partially supported by the NSF grant CMMI-1636193. The second author was partially supported by the National Research Foundation of Ukraine, Grant No.~2020.01/0283.

\bibliographystyle{elsarticle-num}

\end{document}